\documentclass[12pt]{amsart}
\title{Equidimensionality and regularity}
\author{Ying Zong}

\address{Department of Mathematics\\University of Toronto}
\email{zongying@math.utoronto.ca}

\date{} % delete this line to display the current date

\begin{document}

%----------------------------------------------------------------------------------
\maketitle
%\tableofcontents

%\section{}
%\subsection{}

%------------------------------------------------------------------------------------

{\bf Theorem.} --- \emph{Let $S$ be a locally noetherian normal algebraic space of residue characteristics zero with regular formal fibers. Assume that there is a morphism $f:X\to S$, which is locally of finite type, equidimensional and \'{e}tale locally on $S$ admits sections, from a regular algebraic space $X$. Then $S$ is regular and $f$ is flat.}

\begin{proof} The hypothesis on $S$ says that every affine scheme \'{e}tale over $S$ is a noetherian $\mathbf{Q}$-scheme with normal excellent local rings.
\smallskip

The hypothesis on $f$ implies that it is surjective, for \'{e}tale locally on $S$ it admits sections, and that it is universally open (EGA IV 14.4.1), for $S$ is normal. Hence, by EGA O 17.3.3+IV 6.1.5, the statement ``$S$ is regular'' is equivalent to the statement ``$f$ is flat''. 
\smallskip

Suppose given $S'$ a locally noetherian algebraic space with regular formal fibers and given a regular surjective morphism $S'\to S$. Then $S'$ is normal and $X\times_SS'$ is regular, for $S'\to S$ and $X\times_SS'\to X$ are regular morphisms. And $f\times_SS'$ is equidimensional with sections \'{e}tale locally on $S'$. It amounts evidently the same to verifying ``$S'$ is regular and $f\times_SS'$ is flat''. 
\smallskip

In particular, the question is local on $S$ for the \'{e}tale topology. Thus, assume $S$ affine and then local and $f$ given a section. Now $S$ being excellent, its subset $R$ consisting of all those points at which $S$ is regular is open. Localizing at a maximal point of $S-R$ if this set is not empty, one can suppose $S$ of dimension $\geq 2$ and regular outside its closed point $s$. One can suppose furthermore $S$ complete. For, the completion $S'$ of $S$ along $s$ is excellent and the projection $S'\to S$ regular surjective and one replaces $(S, X, f)$ by $(S', X\times_SS', f\times_SS')$.
\smallskip

Observe that $S$ is pure along $s$. Indeed, let $E$ be a locally constant constructible sheaf of sets on $S-\{s\}$ for the \'{e}tale topology. If $j$ denotes the open immersion of $X-f^{-1}(s)$ in $X$, $j_*f^*E$ is locally constant constructible on $X$ (SGA 2 X 3.4), for $X$ is regular and for as $f$ is equidimensional $f^{-1}(s)$ is of codimension $\geq 2$ in $X$. Take $g:S\to X$ a section of $f$ which one assumes is given. Then $g^*j_*f^*E$ extends $E$ to $S$. A similar argument shows that $S$ is also parafactorial along $s$, namely (EGA IV 21.13.11), factorial. (If $S$ is of dimension $2$, purity suffices to imply regularity by Mumford.) 
\smallskip

Let the given section of $f$ be $g: S\to X$. Replacing $X$ by an \'{e}tale neighborhood of $g(s)$, one can suppose $X$ affine and that every irreducible component of $f^{-1}(s)$ contains $g(s)$ and, taking a further finite \'{e}tale base change $S'\to S$ if necessary, that every irreducible component of $f^{-1}(s)$ is geometrically irreducible. 
\smallskip

Choose a maximal point $a$ of $f^{-1}(s)$. Let $\widetilde{X}=\mathrm{Spec}(\mathcal{O}_{X, a})$ and $\overline{X}$ its completion along $a$; the projection $\overline{X}\to \widetilde{X}$ is regular, for $\widetilde{X}$ is excellent. Notice that
\smallskip

--- \emph{The morphism $f: \widetilde{X}-\{a\}\to S-\{s\}$ is regular and bijective }:
\smallskip

 As $S$ is regular outside $s$ and as $X$ is regular, $f$ is smooth along $g(S-\{s\})$. Consider an arbitrary point $v$ of $S-\{s\}$, choose a spectra $V$ of discrete valuation ring and a morphism $V\to S$ such that the generic point $v'$ (resp. closed point $s'$) of $V$ is mapped to $v$ (resp. $s$), and let $X_V=X\times_SV$, $f_V=f\times_SV$, $g_V=g\times_SV$. Since $f_V$ is universally open and is smooth at $g_V(v')$, every sufficiently small open neighborhood $U$ of $g_V(s')$ in $X_V$ has smooth geometrically connected fiber over $v'$; hence, the localization of $U$ at each maximal point of the closed fiber $U_{s'}$ is irreducible of dimension $1$. As $f$ is flat outside $f^{-1}(s)$ and as every irreducible component of $f^{-1}(s)$ is geometrically irreducible, the assertion follows. 
\smallskip

Let $k=k(s)$. By Artin \cite{artin} 4.7 one finds a connected affine $k$-scheme of finite type $T$ and a point $t\in T(k)$ such that $T$ is $k$-smooth outside $t$ and such that $S$ is isomorphic to the completion of $T$ along $t$; by means of this isomorphism $k$ is a chosen coefficient field of $\Gamma(S, \mathcal{O})$. 
\smallskip

As $S$ is normal, $T$ is normal at $t$ and thus normal. As moreover $T$ has a $k$-rational point, it is geometrically connected over $k$. 
Let $k'$ be a coefficient field of the formally smooth $k$-algebra $\Gamma(\overline{X}, \mathcal{O})$. Consider the induced $k'$-algebra homomorphism 
\[\Gamma(T, \mathcal{O})\otimes_kk'\to \Gamma(\overline{X}, \mathcal{O})\] and the $k'$-morphism 
\[p:\overline{X}\to \overline{T'}\] where $\overline{T'}$ is the completion of $T':=T_{k'}$ along $t'$, the unique $k'$-point of $T'$ above $t$; the morphism $\overline{T'}\to S$ induced from $T'\to T$ is regular, for it is formally smooth and $S$ excellent. 
\smallskip

This morphism $p$ is finite by the choice that $a$ is a maximal point of $f^{-1}(s)$, and $p$ is surjective, for $\overline{T'}$ is normal of the same dimension as that of $\overline{X}$. 
\smallskip

If $q: \overline{X}\to \widetilde{X}$ denotes the projection, note that 
\[i:=(q, p): \overline{X}\to \widetilde{X}\times_S\overline{T'}\] is a closed immersion; it is a regular closed immersion since $\overline{X}$ and $\widetilde{X}\times_S\overline{T'}$ are regular. In particular, the conormal module of $i$, $N$ is free over $\mathcal{O}_{\overline{X}}$. On writing $p$ as the composition
\[\overline{X}\stackrel{i}{\longrightarrow}\widetilde{X}\times_S\overline{T'}\stackrel{f\times_S\overline{T'}}{\longrightarrow}\overline{T'}\] one deduces the canonical exact sequence
\[H^{-1}(L_p)\to N\to q^*\Omega^1_{\widetilde{X}/S}\to \Omega^1_p\to 0\] where $L_p$ is the cotangent complex of $p$ and where $H^{-1}(L_p)\to N$ is $0$, the $\mathcal{O}_{\overline{X}}$-module $H^{-1}(L_p)$ (resp. $N$) being torsion (resp. free). By approximation, the injection $N\to q^*\Omega^1_{\widetilde{X}/S}$ is isomorphic to an injection of the form
\[q^*(\alpha): q^*\widetilde{N}\to q^*\Omega^1_{\widetilde{X}/S}\] where $\widetilde{N}$ is a free $\mathcal{O}_{\widetilde{X}}$-module and $\alpha: \widetilde{N}\to \Omega^1_{\widetilde{X}/S}$ an injective $\mathcal{O}_{\widetilde{X}}$-homomorphism. So $\Omega^1_p=q^*\omega$ with $\omega=\mathrm{Coker}(\alpha)$. 
\smallskip

To finish, it is sufficient to show that
\smallskip

--- \emph{The morphism $p: \overline{X}\to \overline{T'}$ is \'{e}tale outside the closed point of $\overline{X}$.}
\smallskip

For then $p$ is \'{e}tale, as the argument proving the purity of $S$ along $s$ applies to $(\overline{T'}, X\times_S\overline{T'}, f\times_S\overline{T'})$, and then, as $\overline{X}$ is regular, $\overline{T'}$ and $S$ are regular.
\smallskip

Write again $a$ (resp. $t'$) for the closed point of $\overline{X}$ (resp. $\overline{T'}$). The finite surjective morphism 
\[p: \overline{X}-\{a\}\to \overline{T'}-\{t'\}\] is flat, since $\overline{T'}-\{t'\}$ and $\overline{X}-\{a\}$ are regular. By SGA 2 X 3.4 it remains to show that $p$ is unramified, that is, $\Omega^1_p=q^*\omega$ vanishes, at every codimension $1$ point of $\overline{X}$. In view that 
\[f: \widetilde{X}-\{a\}\to S-\{s\}\] is regular and bijective, it amounts to showing 
that 
\[V(h\mathcal{O}_{\overline{X}})\cap (\overline{X}-\{a\})\] is reduced for each prime element $h$ of the factorial ring $\Gamma(S, \mathcal{O})$. Both $\overline{X}\to\widetilde{X}$ and $\widetilde{X}-\{a\}\to S-\{s\}$ being regular morphisms, this is clear, hence the proof.

\end{proof}

%-------------------------------------------------------------------------------------------------------------

\bibliographystyle{amsplain}

%-----------------------------------------------------------------

\end{document}